\documentclass[SEDP]{cedram-sem}

\usepackage{graphicx, amssymb}
\usepackage{amsfonts}
\usepackage{amsthm}
\usepackage{xcolor}
\usepackage{hyperref}
\usepackage{bbm}
\theoremstyle{plain}
\newtheorem{theorem}{Theorem}

\theoremstyle{definition}

\title
[]
{Multifractality and polygonal vortex filaments}

\alttitle{Multifractalit\'e et des tourbillons filamentaires polygonaux}

\author
[V. \lastname{Banica}]
{\firstname{Valeria} \lastname{Banica}}

\address{Sorbonne Universit\'e, Universit\'e Paris Cit\'e, CNRS, INRIA, Laboratoire Jacques-Louis Lions, LJLL,\\
F-75005 Paris, France}
\thanks{V.B. is partially
supported by the French ANR grant  BOURGEONS and by the ERC grant GEOEDP}
\email{Valeria.Banica@math.cnrs.fr}

 \author[D. \lastname{Eceizabarrena}]{\firstname{Daniel} \lastname{Eceizabarrena}}
\address{BCAM - Basque Center for Applied Mathematics\\
Bilbao, Spain} 
\thanks{D.E. is partially supported by the European Union Horizon Europe research and innovation programme under Marie Sklodowska Curie Actions with grant agreement 101104250 - TIDE, by the Simons Foundation Collaboration Grant on Wave Turbulence (Nahmod’s award ID 651469), and by the American Mathematical Society and the Simons Foundation under an AMS-Simons Travel Grant for the period 2022-2024.}
\email{deceizabarrena@bcamath.org}

\author[A. R. \lastname{Nahmod}]{\firstname{Andrea R.} \lastname{Nahmod}}
 \address{Department of Mathematics and Statistics\\
 University of Massachusetts Amherst\\
710 N Pleasant St, Amherst, MA 01003, USA}
\thanks{A.N. is partially supported by NSF DMS-2052740, NSF DMS-2101381 and the Simons Foundation Collaboration Grant on Wave Turbulence (Nahmod’s award ID 651469)}
\email{nahmod@umass.edu}

\author[L. \lastname{Vega}]{\firstname{Luis}  \lastname{Vega}}

 \address{Department of Mathematics\\
University of Basque Country\\
Apdo 644, 48080 Bilbao, Spain}

\thanks{L.V. is partially supported by MICINN (Spain) CEX2021-001142, PID2021-126813NB-I00 (ERDF A way of making Europe) and IT1247-19 (Gobierno Vasco).}
\email{lvega@bcamath.org}

\keywords{Vortex filaments, multifractality, Riemann's non-differentiable function, Diophantine approximation.}
  
\altkeywords{Exemple, Annales, classe \LaTeX}

\subjclass{11J82, 11J83, 26A27, 28A78, 42A16, 76F99}

\begin{document}


\begin{abstract}
 In this proceedings article we survey the results in \cite{BanicaEceizabarrenaNahmodVega2024} and their motivation, as presented at the 50th \currentjournaltitle \,2024. With the aim of quantifying turbulent behaviors of vortex filaments, we study the multifractality of a family of generalized Riemann's non-differentiable functions. These functions represent, in a certain limit, the trajectory of regular polygonal vortex filaments that evolve according to the binormal flow, the classical model for vortex filaments dynamics. We explain how we determined their spectrum of singularities
through a careful design of Diophantine sets, which we study by using the Duffin-Schaeffer theorem and the Mass Transference Principle.\end{abstract}

\begin{altabstract}
Dans cet acte de conf\'erence nous passons en revue les r\'esultats de \cite{BanicaEceizabarrenaNahmodVega2024} et leur motivation, tels qu'ils ont \'et\'e pr\'esent\'es au 50e \currentjournaltitle \,2024. Dans le but de quantifier les comportements turbulents des filaments tourbillonnaires, nous \'etudions la multifractalit\'e d'une famille de fonctions non diff\'erentiables de Riemann g\'en\'eralis\'ees. Ces fonctions repr\'esentent, dans une certaine limite, la trajectoire de filaments tourbillonaires polygonaux r\'eguliers qui \'evoluent selon le flot binormal, le mod\`ele classique pour la dynamique des tourbillons filamentaires. Nous expliquons comment nous avons d\'etermin\'e pour certaines de ces fonctions le spectre des singularit\'es. La preuve repose sur une construction d'ensembles diophantiens que nous \'etudions en utilisant le th\'eor\`eme de Duffin-Schaeffer et le principe de transfert de masse.

\end{altabstract}

\maketitle

%
%


\section{Multifractal analysis and Riemann's non-differentiable function} 

Multifractal analysis started in the 1980s to explain the deviations from Kolmogorov's 1941 theory of turbulence observed by Anselmet et al. \cite{AnselmetGagneHopfingerAntonia1984} (see Arneodo and Jaffard's expository article \cite{ArneodoJaffard2004}). 
In these experiments, 
the speed of the turbulent fluid seems very irregular in some regions and less in others, and the borders between regions are unclear.
What is more, zooming in very irregular regions one finds less irregular regions, and vice-versa. 
Facing such a non-uniform distribution of regularities, instead of explicitly computing the regularity at each individual point it turns out that it is more suitable to measure the sets where a given regularity is reached.
A typical way to formalize this is as follows.

For $f:\mathbb R\rightarrow\mathbb R$ and $0<\alpha<1$, the local H\"older regularity of $f$ at a point $t$ is measured as
\[
f\in \mathcal C^\alpha(t) \quad \iff \quad |f(t+h)-f(t)|\leq Ch^\alpha,\quad  \forall h\ll 1,
\]
and the local H\"older exponent of $f$ at a point $t$ is given by
\[
\alpha_f(t)=\sup\{ \, \alpha \, : \,  f\in \mathcal C^\alpha(t) \, \}.
\]
A  multifractal function is a function possessing infinitely many local H\"older exponents. The multifractal analysis of a function consists in computing its spectrum of singularities,
that is, 
the Hausdorff dimension\footnote{We recall that the $s$-Hausdorff measure of a set is
\begin{equation*}
\mathcal {H}^s(A)=\lim_{\delta\rightarrow 0} \left( \inf \Big\{ \sum_{j=1}^\infty \operatorname{diam} (A_j)^s \, : \, A\subset \bigcup_{j=1}^\infty A_j, \, \,  \operatorname{diam}(A_j) < \delta \,  \Big\} \right), 
\end{equation*}    
and the Hausdorff dimension of a set $A$ is 
$$\dim_{\mathcal H}A=\sup \{ s \, : \, \mathcal{H}^s(A)=\infty \, \}=\inf\{ \, s\, : \, \mathcal{H}^s(A)=0\}.$$} 
of its iso-H\"older sets
$$d_f(\alpha)=\dim_{\mathcal H}\{ \,  t \, : \, \alpha_f(t)=\alpha \, \}.$$
By convention,  we set $d_f(\alpha)=-\infty$ if the set $\{ \,  t \, : \, \alpha_f(t)=\alpha \, \}$ is empty. 

A classical example of multifractal function is
$$R(t)=\sum_{n =1}^\infty \frac{\sin(n^2 t)}{n^2}, $$
proposed by Riemann in 1860 when he was looking for continuous functions that are nowhere differentiable. 
Weierstrass, who failed to prove the claim that $R$ is actually nowhere differentiable (which eventually turned out to be false!), 
proved that the function
$$W(t)=\sum_{n=1}^\infty \frac{\cos(4^n t)}{2^n}, $$
which is simpler to study, is nowhere differentiable. 
However, 
it is not multifractal because it has H\"older exponent $1/2$ everywhere; it is rather a \textit{monofractal} function. 
In particular, its spectrum of singularity is
\begin{equation*}
d_W(\alpha) = 
\left\{
\begin{array}{ll}
1, & \alpha=1/2,\\
-\infty, &  \alpha\neq 1/2.
\end{array} 
\right.
\end{equation*}
Riemann's function is richer; 
Hardy \cite{Hardy1916} and Gerver \cite{Gerver1970,Gerver1971}
showed that it actually allows a derivative at certain points, 
and it was proved by Jaffard \cite{Jaffard1996} to be a genuine multifractal function 
with spectrum of singularities
\begin{equation*}
d_R(\alpha)=
\left\{
\begin{array}{ll}
4\alpha-2, & 1/2 \leq \alpha \leq 3/4, \\ 
0, &  \alpha=3/2,\\
-\infty, &  \text{otherwise,}
\end{array}
\right.
\end{equation*}
a result that holds too for its complex version given by the function
$$R_0(t)=\sum_{n \in \mathbb Z \setminus \{0 \} }\frac{e^{in^2 t}}{n^2}.$$
Jaffard actually computed the local H\"older regularity at every point, 
a result for which Broucke and Vindas \cite{BrouckeVindas2023} gave recently an alternative proof.

\begin{figure}
\includegraphics[width=0.32\linewidth]{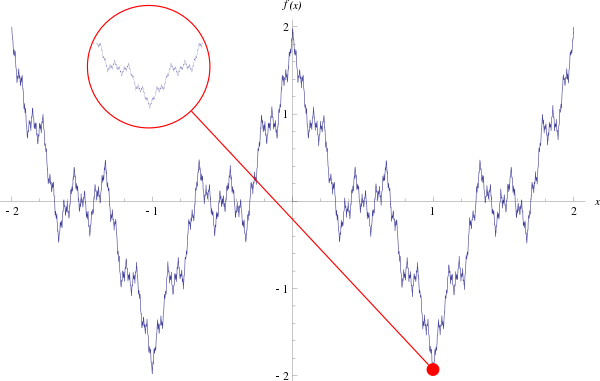}
\includegraphics[width=0.32\linewidth]{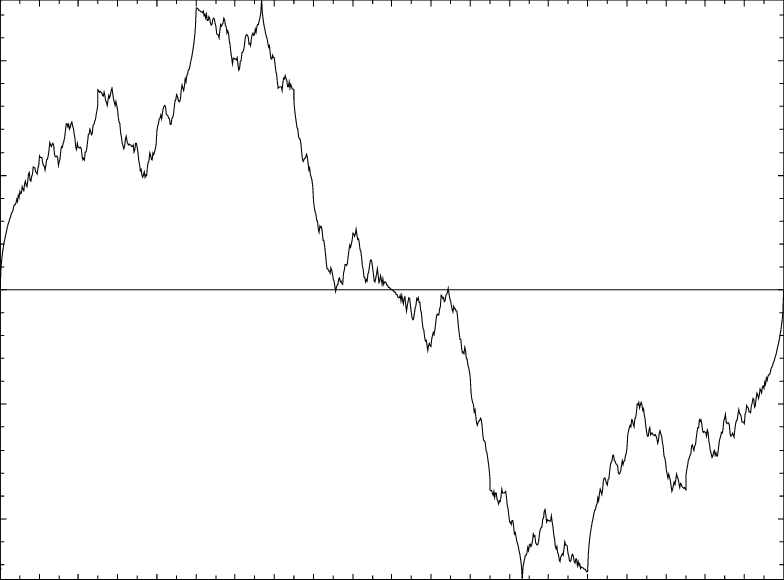}
\includegraphics[width=0.32\linewidth]{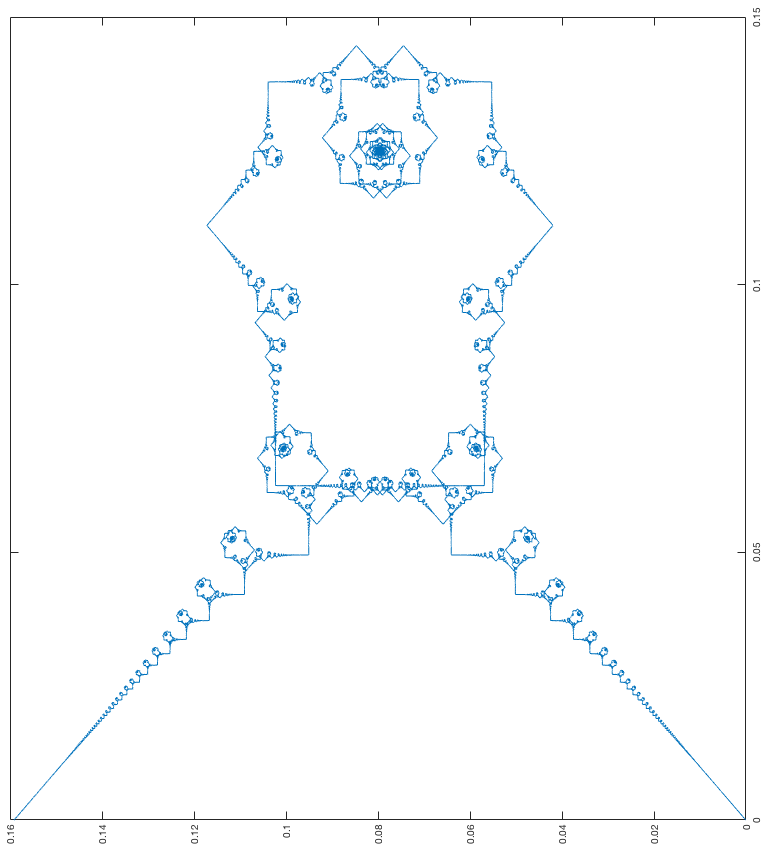}
\caption{The graphs of $W$, $R$, and the image of $\widetilde R_0$.}
\end{figure}

Other notions related to turbulence were successfully tested on Riemann's function. Jaffard \cite{Jaffard1996} showed that it satisfies the Frisch-Parisi multifractal formalism (see \eqref{FP}), and Boritchev, Vila\c{c}a Da Rocha and Eceizabarrena \cite{BoritchevDaRochaEceizabarrena2021} proved that it is intermittent (see \eqref{flat}). 
Also, Eceizabarrena \cite{Eceizabarrena2020,Eceizabarrena2021} showed that the image of $R_0$ has no tangents and has Hausdorff dimension\footnote{Note that even the Hausdorff dimension of the graph of the Weierstrass function was determined only recently by Shen \cite{Weixiao2018}} less than $4/3$. 

Generalizations of Riemann's function were intensively studied after the work of Jaffard. 
For example, 
Chamizo and Ubis \cite{ChamizoUbis2014}, and Seuret and Ubis \cite{SeuretUbis2017} obtained upper and lower bounds for the spectrum of singularities of $\sum_n\frac{e^{iP(n)t}}{n^\alpha}$ with $P$ a polynomial, 
with an approach different from Jaffard's wavelet techniques and the unimodular group action on Riemann's function.
Vila\c{c}a Da Rocha and Eceizabarrena \cite{DaRochaEceizabarrena2022} studied the intermittency of $\sum_n \frac{e^{in^2t}}{n^\alpha}$ with $\alpha > 1/2$. 
Kapitanski and Rodnianski \cite{KapitanskiRodnianski1999} gave fine regularity results for $\sum_n e^{in^2t_0+inx}$ with $t_0$ fixed, 
which represents the fundamental solution to the periodic Schr\"odinger equation at time $t_0$. 
Banica and Vega \cite{BanicaVega2022} proved that the spectrum of singularities of 
$\sum_{n \in \mathbb N} \frac{e^{i(pn+q)^2t}}{(pn+q)^2}$ for $p,q\in\mathbb N$ is the same as for  Riemann's function, and that the Frisch-Parisi multifractal formalism is also satisfied. 

\section{Main result}

\begin{figure}
\includegraphics[width=0.7\linewidth]{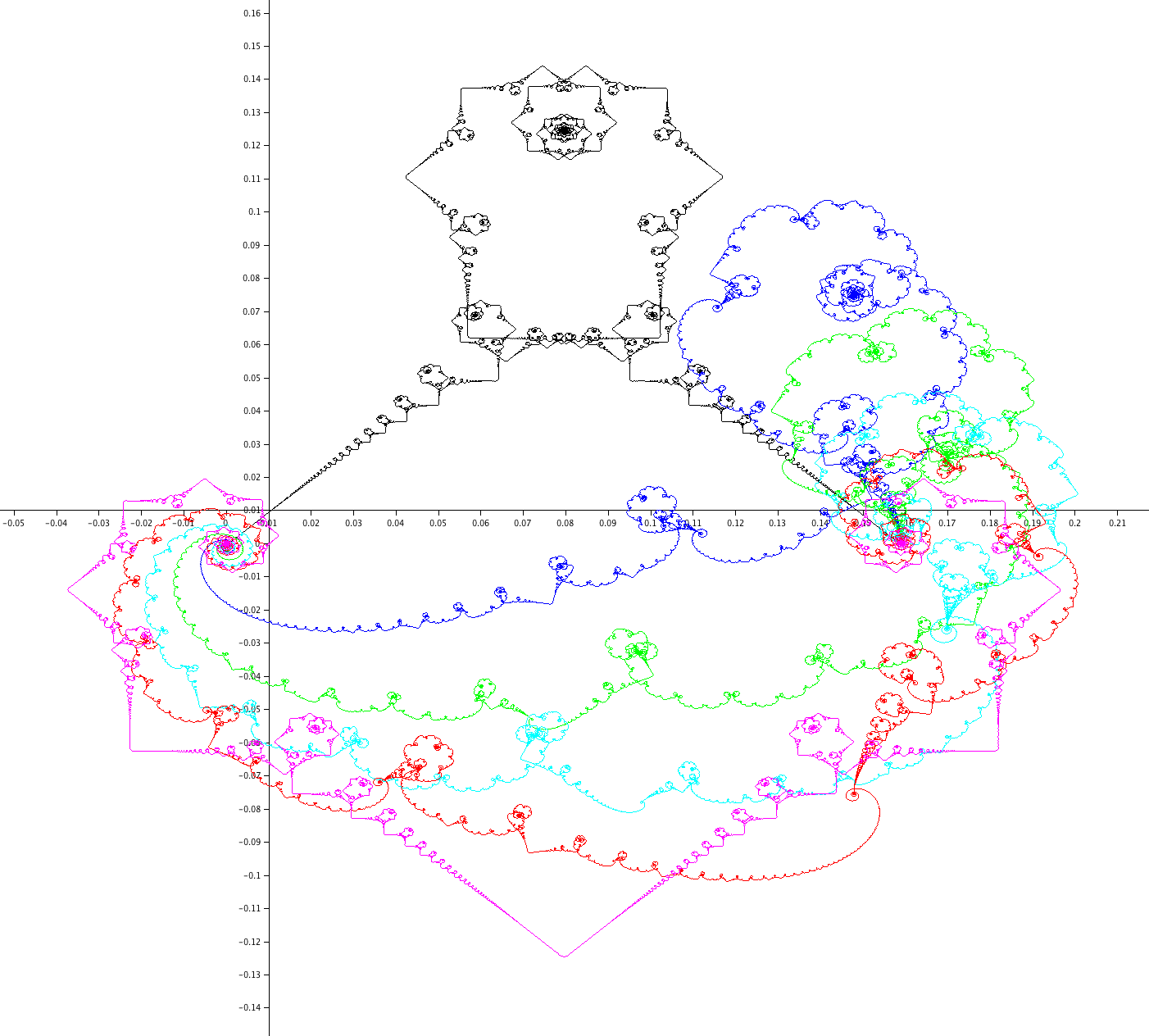}
\caption{The curve $\widetilde R_{x_0}([0,2\pi])$ for several values of $x_0$: $0$ (black), $0.1$ (blue), $0.2$ (green), $0.3$ (cyan), $0.4$ (red), $0.5$ (magenta).}
\label{fig:Trajectories}
\end{figure}

In \cite{BanicaEceizabarrenaNahmodVega2024} we study a different generalization of the Riemann function, namely 
$$R_{x_0}(t)=\sum_{n \in  \mathbb Z \setminus \{ 0 \}  } \frac{e^{2 \pi i (n^2t+nx_0)}}{n^2}, \qquad \text{ for } x_0 \text{ fixed. } $$
These functions, or rather the closely related\footnote{In particular, the functions $R_{x_0}$ and $\widetilde{R}_{x_0}$ have the same regularity in $t$.}   
\[
\widetilde{R}_{x_0}(t) = \sum_{n \in  \mathbb Z  } \frac{e^{in^2t} - 1}{n^2} \, e^{i n x_0},
\] 
arise naturally in the setting of vortex filaments, as we explain in Section~\ref{SectionVFE}.
In Figure~\ref{fig:Trajectories} we display the images of $\widetilde{R}_{x_0}$ for different values of $x_0$. 
The result from \cite{BanicaEceizabarrenaNahmodVega2024} that we here review is the following.
 
\begin{theorem}\label{th}
Let $x_0\in\mathbb R$. 
Then, the function $R_{x_0}$
is multifractal, with infinitely many local H\"older exponents. 
If $x_0\in\mathbb Q$, 
the spectrum of singularities of $R_{x_0}$ is
\begin{equation}
d_{R_{x_0}}(\alpha) = 
\left\{\begin{array}{ll}
4\alpha-2, &  1/2 \leq \alpha \leq 3/4,\\
 0, &  \alpha= 3/2, \\
 -\infty, & \text{otherwise.}
 \end{array}\right.
\end{equation}
\end{theorem}
  
The proof starts as in Jaffard's \cite{Jaffard1996} and then it follows the approach by Chamizo and Ubis in \cite{ChamizoUbis2014}, 
but it ends up with new Diophantine sets that approximate the iso-H\"older sets. We measure these new sets combining the Duffin-Schaeffer theorem\footnote{
The Duffin-Schaeffer theorem was conjectured by Duffin and Schaeffer \cite{DuffinSchaeffer1941} in 1941,
and they proved it in some particular cases that suffice to prove Theorem~\ref{th} for $x_0 \in \mathbb Q$. It was recently proved in all generality by Koukoulopoulos and Maynard \cite{KoukoulopoulosMaynard2020}, 
which we need for our results in \cite{BanicaEceizabarrenaNahmodVega2024} when $x_0 \not\in \mathbb Q$.
}
\cite{DuffinSchaeffer1941}
and the Beresnevich-Velani Mass Transference Principle \cite{BeresnevichVelani2006}. 
We explain this procedure in detail in Section~\ref{SectionProof}. 

In \cite{BanicaEceizabarrenaNahmodVega2024} we also showed that $R_{x_0}$ is intermittent in small scales by proving that its flatness tends to infinity when the scale parameter tends to zero, that is, 
\begin{equation}\label{flat}
F_{R_{x_0}}(N):=\frac{\|P_{\geq N}R_{x_0}\|_4^4}{\|P_{\geq N}R_{x_0}\|_2^4}\overset{N\rightarrow\infty}{\longrightarrow}\infty,
\end{equation}
where $P_{\geq N}$ is the high-pass filter of Fourier modes larger than $N$. 
As a consequence, one can deduce that the function satisfies the Frisch-Parisi multifractal formalism\footnote{
The Frisch-Parisi multifractal formalism was originally proposed for the velocity in an Eulerian setting, 
but it can be equally proposed in the Lagrangian setting, 
which Riemann's function is in principle more related to since it represents a time trajectory.
 See the work of Chevillard et al. \cite{ChevillardCastaignArneodoLevequePintonRoux2012} for a discussion on the differences between these two frameworks.}
\begin{equation}\label{FP}
d_{R_{x_0}}(\alpha) = \inf_p \left\{  \alpha p-\eta_{R_{x_0}}(p)+1 \right\},
\qquad \text{ for } \quad  \frac12 \leq \alpha \leq \frac34,
\end{equation}
where $\eta_{R_{x_0}}(p):=\sup\{s,\,R_{x_0}\in B^{ s/p}_{p,\infty}\}$
and $B^{ s/p}_{p,\infty}$ stands for the Besov space\footnote{See the recent works by Barral and Seuret \cite{BarralSeuret2023,BarralSeuret2023bis} on the validity of the multifractal formalism in Besov spaces.}. 

Regarding the case when $x_0$ is irrational, 
we proved that $R_{x_0}$ is multifractal by giving a result on its spectrum of singularities, although not recovering it explicitly. 
The difficulty comes from the interference between the exponents of irrationality of both $x_0$ and $t$.

\section{Motivation from fluid dynamics}\label{SectionVFE}

In this section we explain how Riemann's function $R_{x_0}$ appears naturally in the evolution of curves that evolve according to the binormal flow model for vortex filaments dynamics. 
Let us first recall this model and some very recent results related to it. 

\subsection{The binormal flow}
The binormal flow (BF), also known as local induction approximation (LIA) or vortex filament equation (VFE), 
is the oldest, simplest and richest formally derived model for one vortex filament dynamics in a 3D fluid governed by Euler's equations. If the vorticity at time $t$ is concentrated along an arclength-parametrized curve $\chi(t)$ in $\mathbb R^3$, its evolution in time is expected to evolve according to the equation
\begin{equation}
\label{BF}
\chi_t=\chi_x\times\chi_{xx}, \qquad \text{ or equivalently } \qquad\chi_t  =c \,b,
\end{equation}
where $b$ is the binormal vector of the curve $\chi$ and $c$ is its curvature. 
This model was derived formally from the Biot-Savart formula by Da Rios in 1906 \cite{DaRios1906}, 
and it was justified rigorously by Jerrard and Seis in 2017 \cite{JerrardSeis2017}.
Understanding if and when the vorticity propagates its initial structure of being concentrated along a curve is still a very difficult open problem. 

Let us briefly recall a few very recent advances in this direction.
Concerning the Navier-Stokes equation, 
Bedrossian, Germain and Harrop-Griffiths \cite{BedrossianGermainHarropGriffiths2023} proved that the Cauchy problem is locally well-posed for an initial filament data  with no symmetry assumptions, 
but for times that are too small to observe the binormal flow and to pass to the vanishing viscosity limit. 
The vanishing viscosity limit was proved by Gallay and Sverak \cite{GallaySverak2024} for the particular case of axisymmetric vortex rings. 
The BF dynamics was recovered by Fontelos and Vega \cite{FontelosVega2023} for Giga-Miyakawa solutions with initial filament data with no symmetry assumptions, but this regime does not allow to pass to the vanishing viscosity limit. 
For the Euler equations Donati, Lacave and Miot \cite{DonatiLacaveMiot2024}, and previously D\'avila, Del Pino, Musso and Wei \cite{DavilaDelPinoMussoWei2022} with other methods, constructed solutions with vorticity concentrated on helices, which are particular solutions of the binormal flow. These are configurations with symmetries benefitting from a dimensional reduction. 
Thus, despite recent efforts, the binormal flow conjecture is still a serious gap away from being understood.

\subsection{Experiments, numerics, and a rigorous result}

\begin{figure}[h]
\includegraphics[width=0.49\linewidth]{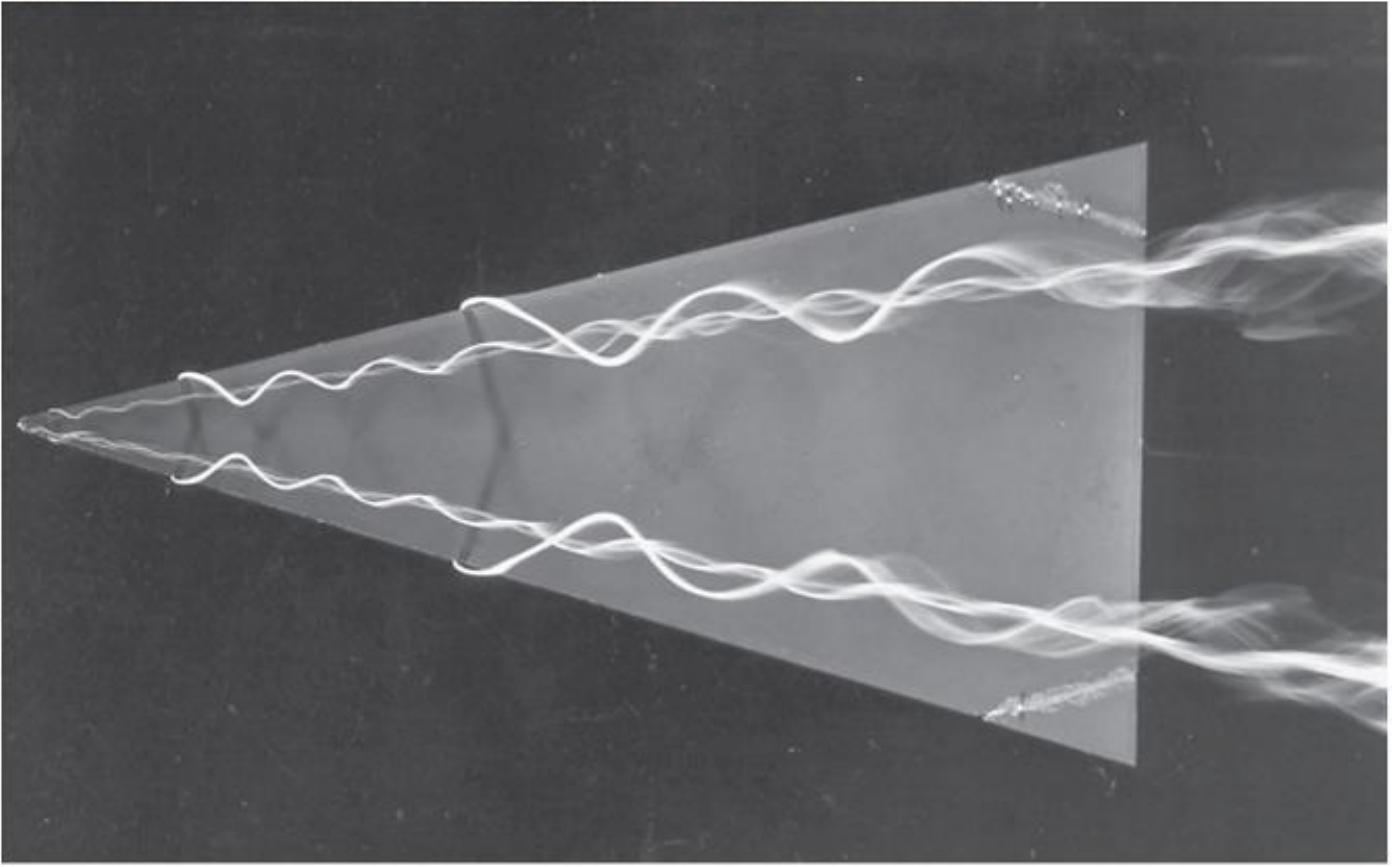}
\includegraphics[width=0.368\linewidth]{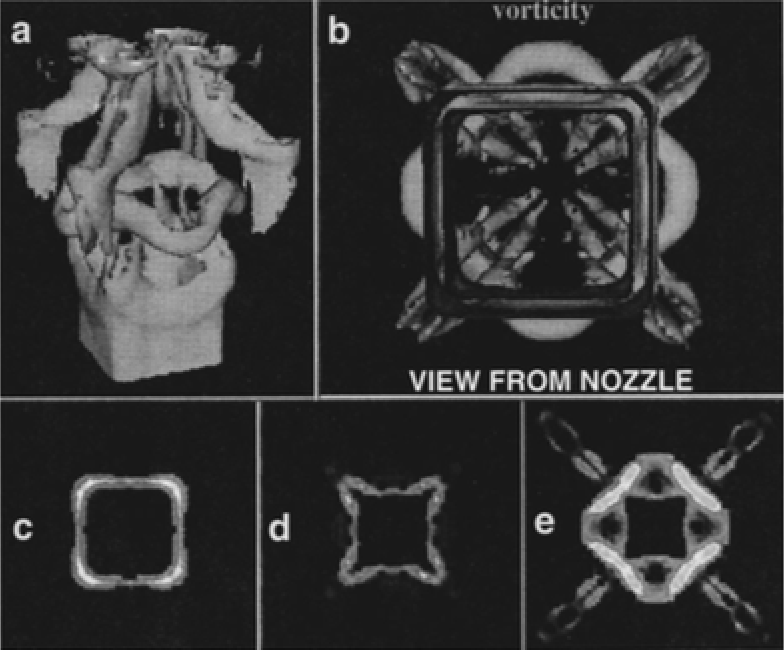}
\caption{On the left: vortices in a fluid flowing over a triangular obstacle. On the right: axis switching in numerical simulation of square jets.}
\label{fig:Vortices_Axix_Switching}
\end{figure}

A special class of solutions of the binormal flow are the self-similar solutions, which are smooth curves that develop a singularity in the shape of a corner  in finite time. They were known and used by physicists since the 1980s in the framework of reconnection of vortex filaments in ferromagnetics, but they were not rigorously studied until 2003 by Guti\'errez, Rivas and Vega \cite{GutierrezRivasVega2003}. This type of dynamics appears in fluids passing over a triangular obstacle and in trail vortex reconnection, as shown in Figure~\ref{fig:Vortices_Axix_Switching} left. 
The interaction of many corner singularities yields a range of complex behaviors such as the Talbot effect, $L^\infty$-cascades of energy, rogue waves, intermittency and multifractal behaviors \cite{BanicaEceizabarrenaNahmodVega2024,BanicaVega2020,BanicaVega2022,BanicaVega2022bis,BanicaVega2024}. Here we consider the last type of results that are the object of this proceedings article. 

\begin{figure}[h]
\includegraphics[width=0.8\linewidth]{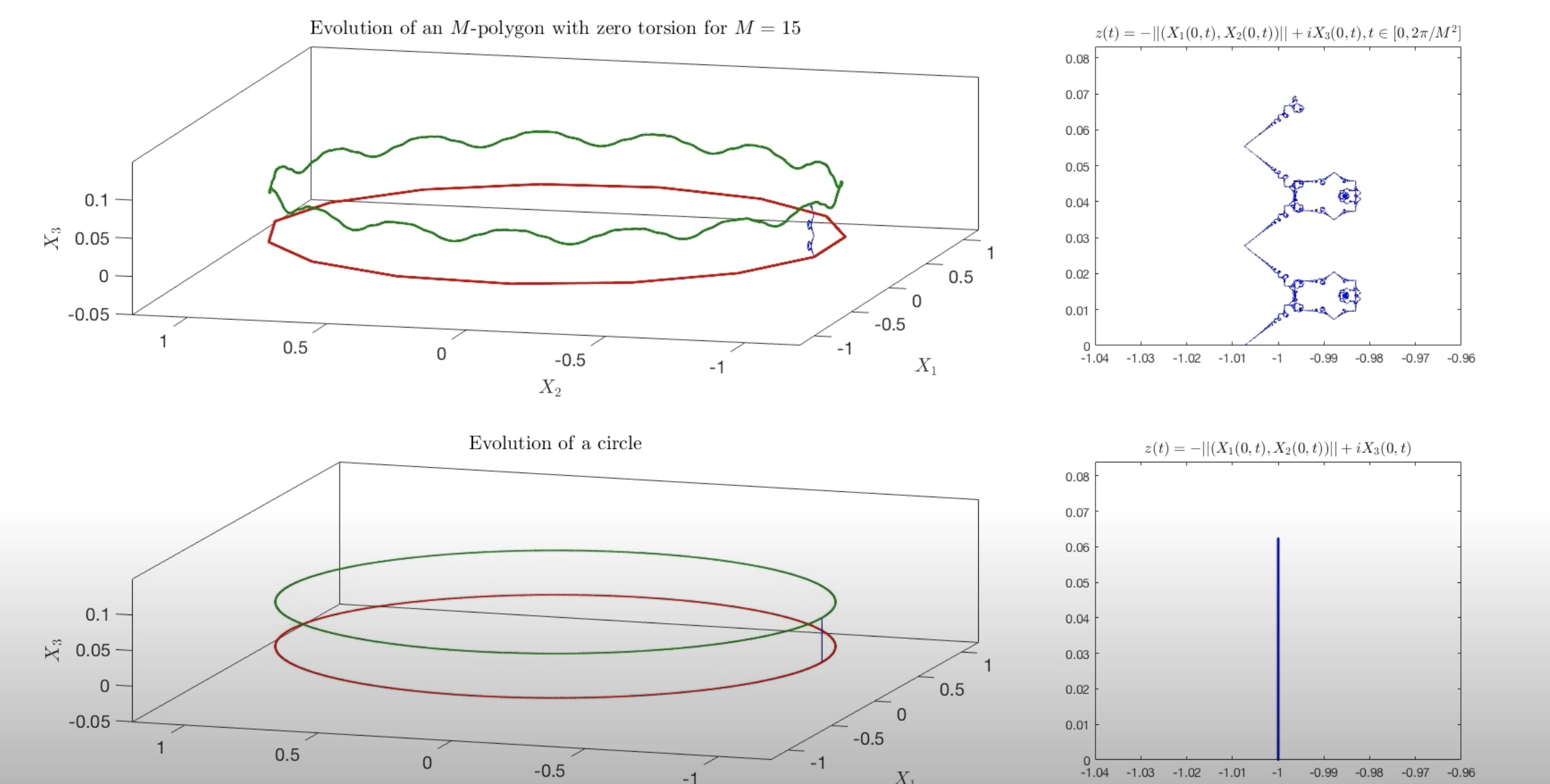}
\caption{On the left: Evolution by the binormal flow of a circle and of an $M$-polygon with $M=15$. On the right: the trajectory in time of the solution at $x=0$.  }
\label{fig:Video}
\end{figure}

Noncircular jets such as square jets have been studied since the 1980s for the turbulent features they produce. 
The experiments by Todoya and Hussain \cite{TodoyaHussain1989}
and the numerics by Grinstein and De Vore \cite{GrinsteinDeVore1996} are some examples of this, see Figure 3 right. 
At the level of the binormal flow this corresponds to considering as initial data a closed curve with the shape of a regular polygon. 
Such a regular polygon of $M$ sides $\chi_M(0,x)$ with corners located\footnote{We see here the closed curve as being parametrized, in a periodic way, by $x\in\mathbb R$} at $x\in\mathbb Z$ is expected to evolve by the binormal flow to skew polygons of $Mq$ sides at times $p/q \in\mathbb Q$, 
as suggested by numerics by Jerrard and Smets \cite{JerrardSmets2015} and by De la Hoz and Vega \cite{delaHozVega2014}. 
Moreover, the trajectories of the corners $\chi_M(t,0)$ were numerically proved to behave like Riemann's function $R_{0}(t)$ when $M\rightarrow\infty$ by De la Hoz and Vega \cite{delaHozVega2014} (see also De la Hoz, Kumar and Vega \cite{delaHozKumarVega2020} and the corresponding video \cite{delaHozKumarVegaVideo}, a screenshot of which we show in Figure~\ref{fig:Video}). 
The presence of Riemann's function in this context was rigorously proved by Banica and Vega \cite{BanicaVega2022} by taking as initial data polygonal line approximations made by regular $M$-polygons (looped a large amount of times) and two half-lines, 
for which they used the solutions constructed in their previous work \cite{BanicaVega2020}
This approach extends to the trajectories of all locations $x_0$, giving rise to $R_{x_0}$.

\subsection{Appearance of Riemann's function}\label{Sec_Appearance}

Let us explain with more details how $R_{x_0}$ arises in this context. 
To do so, let us first briefly state the links between the binormal flow and Schr\"odinger equations. 
One can immediately see by taking the derivative in $x$ in \eqref{BF} that the tangent vector $T$ of a smooth solution $\chi$ to the binormal flow satisfies the Schr\"odinger map equation from $\mathbb R$ to $\mathbb S^2$, which is the Heisenberg ferromagnetic continuous model
\begin{equation}\label{SM}
T_t=T\wedge T_{xx}.
\end{equation}
On the other hand, in 1972 Hasimoto \cite{Hasimoto1972} introduced a transform\footnote{This transform can be seen as an inverse Madelung transform.} that connects explicitly the binormal flow and the 1D cubic Schr\"odinger equation (NLS). More precisely, for curves with non-vanishing curvature, the function
$$u(t,x)=c(t,x)e^{i\int_0^x \tau(t,y)dy},$$
where $c$ and $\tau$ are respectively the curvature and torsion of a smooth solution to the BF, solves
\begin{equation}\label{NLS}
i\partial_tu+\Delta u+\frac 12(|u|^2-A(t))u=0,
\end{equation}
where $A(t)$ is a function that depends only on time and is defined in terms of the curvature and torsion. 
Koiso  \cite{Koiso1997} showed in 1997 that the non-vanishing condition on the curvature can be removed by using parallel transport frames instead of Frenet frames, 
that is, by constructing an orthonormal basis $(T,e_1,e_2)$ verifying
\begin{equation}
\label{PT}
\left(  
\begin{array}{c}
T \\
e_1 \\
e_2
\end{array}
\right)_x
= 
\left(  
\begin{array}{ccc}
0 & \alpha & \beta  \\
-\alpha & 0 & 0 \\
-\beta & 0 & 0
\end{array}
\right)
\, 
\left(  
\begin{array}{c}
T \\
e_1 \\
e_2
\end{array}
\right)
\end{equation}
and 
\begin{equation}
\label{PTtgen}
\left(  
\begin{array}{c}
T \\
e_1 \\
e_2
\end{array}
\right)_t
= 
\left(  
\begin{array}{ccc}
0 & -\beta_x & \alpha_x  \\
\beta_x & 0  & -\frac12 (\alpha^2 + \beta^2 - A(t))  \\
-\alpha_x & \frac12 (\alpha^2 + \beta^2 - A(t)) & 0
\end{array}
\right)
\, 
\left(  
\begin{array}{c}
T \\
e_1 \\
e_2
\end{array}
\right)
\end{equation}
In this case, $u=\alpha + i \beta$ solves \eqref{NLS}, possibly with a different $A(t)$.

Conversely, from a smooth solution $u$ to NLS in \eqref{NLS} with any\footnote{This general choice is possible due to a gauge invariance.} $A(t)$,
 one can construct a smooth solution to the BF with tangent vector $T$ satisfying the Schr\"odinger map \eqref{SM}, and frames $(T,e_1,e_2)$ satisfying \eqref{PT}-\eqref{PTtgen} with $\alpha= \operatorname{Re}(u)$ and $\beta= \operatorname{Im}(u)$. 
Thus, one could generate singularities for BF by generating singularities for smooth NLS solutions.
For instance, Guti\'errez, Rivas and Vega \cite{GutierrezRivasVega2003} showed that the self-similar solutions of the BF are obtained from solutions of the type $t^{-1/2} \, e^{i \, x^2/4t} =e^{it\Delta}\delta$ of \eqref{NLS} with $A(t) = 1/t$. 
When this solution, which is smooth for $t>0$, generates the Dirac-singularity at $t=0$, the associated smooth BF solution generates a corner. 
The BF evolutions $\chi_M$ of polygonal-line approximations of $M$-regular polygons were constructed by Banica and Vega \cite{BanicaVega2020} from NLS evolutions of truncations of the Dirac comb $\sum_{|n|\lesssim M} \delta_n$. These NLS evolutions  were proved to be truncations of type $\Sigma_{|n|\lesssim M}e^{it\Delta}\delta_n$ plus a remainder term. By using  the Poison summation formula they rewrite as truncations of type $\Sigma_{|n|\lesssim M}e^{in^2t+inx}$ plus a remainder term. 
In view of \eqref{PT}, 
the derivative in time  $\partial_t\chi_M=T_M \wedge \partial_x T_M$ involves the corresponding NLS solution.  
Consequently, the leading term of the trajectories $\chi_M(t,x_0)$ is proved to be 
$$\int_0^t\sum_{n\in\mathbb Z} e^{in^2\tau+inx_0} \, d\tau=-i \widetilde R_{x_0}(t),$$
whence Riemann's function naturally arises.

\subsection{Energy cascades}
Let us briefly mention another consequence of the construction in \cite{BanicaVega2020}:
the existence of an energy cascade for the tangent vector $T$, 
which can also be understood as an energy cascade for a linear Hamiltonian system. 
This was shown in \cite{BanicaVega2022bis} 
using the method in Section~\ref{Sec_Appearance} 
and solutions of \eqref{NLS} with $A(t)\equiv 0$ of the type  
\[
u(t,x)=\frac{1}{(it)^{1/2}}e^{i\frac{|x|^2}{4t}}\overline V\left(\frac 1t, \frac{x}{t}\right),
\]
with $V(x,t)$ a smooth function that is periodic in $x$. 
Therefore, $(T,e_1,e_2)$ solves \eqref{PT} which, due to $u = \alpha + i\beta$, becomes
\begin{equation}
\label{PTt}
\left(  
\begin{array}{c}
T \\
e_1 \\
e_2
\end{array}
\right)_t
= 
\left(  
\begin{array}{ccc}
0 & -\operatorname{Im} u_x & \operatorname{Re} u_x  \\
\operatorname{Im} u_x & 0  & -\frac12 |u|^2 \\
-\operatorname{Re} u_x  & \frac12 |u|^2 & 0
\end{array}
\right)
\, 
\left(  
\begin{array}{c}
T \\
e_1 \\
e_2
\end{array}
\right).
\end{equation}
The matrix coefficients involve 
\[|u(t,x)|^2= \frac{1}{t}\Big| V\Big(\frac{1}{t}, \frac{x}{t} \Big)\Big|^2 \]
which is a real potential, and also the derivatives in $x$ given by
\begin{equation}\label{Transport_Term}
u_x(t,x)= i\frac{x}{2t}u(t,x)  + \text{Remainder term}.
\end{equation}  
The first term comes from the space derivative of the exponential, and the remainder term corresponds to the space derivative hitting the smooth periodic in space function $V$. Therefore, with this particular choice of $u$, the equation for $T_t$ in \eqref{PTt} can be seen as a linear equation with variable coefficients whose leading behaviour is given by the potential $i\frac{x}{2t}u$. The multiplication by $i \frac{x}{2t}$ creates a travelling wave in Fourier space.
We refer to Theorem 1.1 in \cite{BanicaVega2022bis} for the details.

This phenomenon is reminiscent of the works of 
Apolin\'ario et al. \cite{ApolinarioBeckChevillardGallagherGrande2023,ApolinarioChevillardMourrat2022}
where they proposed abstract linear equations that mimic the phenomenology of energy cascades when the external force is a statistically homogeneous and stationary stochastic process. 
Indeed, these equations have a potential of the type $i c x$  with $c \in \mathbb R$ (see (2.2) in \cite{ApolinarioBeckChevillardGallagherGrande2023}), 
independent of time.
Existence of energy cascades for linear systems have been also proved by Colin de Verdi\`ere and Saint-Raymond in \cite{deVerdiereSaintRaymond2020}.




\section{Sketch of the proof of Theorem \ref{th}}\label{SectionProof}
In this last section we give an overview of the proof of Theorem \ref{th} for $\alpha\in[1/2, 3/4]$. 
Further details can be found in \cite{BanicaEceizabarrenaNahmodVega2024}.   
To determine the spectrum of singularities of $R_{x_0}$ we first estimate variations near rational $t$. 
This will allow to understand the variations near irrationals, depending on how well they can be approximated by rationals. 
By doing so, we will capture every iso-H\"older set between a larger set $\mathbf A$ and a strictly smaller subset of irrationals $\mathbf B$ defined by a constrained Diophantine condition that depends on $x_0$. 
The Hausdorff dimension of the larger set $\mathbf A$ is settled directly by the Jarn\'{i}k-Besicovitch theorem. 
The Hausdorff dimension of the smaller set $\mathbf B$ cannot be obtained directly with the arguments that Jaffard \cite{Jaffard1996} and Chamizo and Ubis \cite{ChamizoUbis2014} used for $x_0=0$, where the Diophantine condition is just a parity condition on the denominators. 
Instead, we will use first the Duffin-Schaeffer theorem \cite{DuffinSchaeffer1941,KoukoulopoulosMaynard2020} to compute the Lebesgue measure of the sets $\mathbf B$.
Once we know that, since these sets are limsups of balls, 
we will use the Mass Transference Principle \cite{BeresnevichVelani2006} to give a lower bound for their Hausdorff dimension\footnote{Which will depend on how much the balls must be dilated so that the limsup of such dilated balls has Lebesgue measure 1.}.

\subsection{Variation of $R_{x_0}$ at rationals}
To compute the behavior of $R_{x_0}$ around a rational $p/q\in\mathbb Q$, 
add the $n=0$ term in the sum, 
split it modulo $q$, and use the Poisson summation formula to get
\begin{equation*}
\begin{split}
R_{x_0}\Big(\frac pq+h\Big)-R_{x_0}\Big(\frac pq\Big)
& = \sum_{n \in\mathbb Z} e^{2 \pi i n^2\frac pq }\frac{e^{2\pi in^2h}-1}{n^2}e^{2\pi inx_0} - 2\pi ih \\
& = \frac{\sqrt{h}}{q}\sum_{m\in\mathbb Z}G(p,m_{x_0,q}+m,q) F\left(\frac{ \operatorname{dist}(x_0,\frac{\mathbb Z}{q})-\frac mq}{\sqrt{h}}\right) - 2\pi i h,
\end{split}
\end{equation*}
where $\operatorname{dist}(x_0,\frac{\mathbb Z}{q}) = x_0 - \frac{m_{x_0,q}}{q}$ and
$F(x)=\mathcal F(\frac{e^{2\pi i\xi^2}-1}{\xi^2})(x) = O(\frac 1{x^2})$
with 
 $F(0)\neq 0$. Also, 
$$G(p,b,q)=\sum_{r=0}^{q-1}e^{2\pi i\frac{p\, r^2 + b\, r}{q}}$$ 
are Gauss sums whose modulus is $\sqrt{q}$ 
except when $q$ is even and
$q/2$ and $b$ have different parity, in which case it is zero. 
Due to the decay of $F$, the leading term is $m=0$, so
$$
R_{x_0}\Big(\frac pq+h\Big)-R_{x_0}\Big(\frac pq\Big)=\frac{\sqrt{h}}{q}G(p,m_{x_0,q},q) F\left(\frac{\operatorname{dist}(x_0,\frac{\mathbb Z}{q})}{\sqrt{h}}\right) - 2\pi ih+O(\min\{\sqrt{q}h,q^\frac 32 h^\frac 32\}).
$$
From here, we obtain a general upper bound valid for all $x_0$ given by
\begin{equation}\label{varrat}
\Big|R_{x_0}\Big(\frac pq+h\Big)-R_{x_0}\Big(\frac pq\Big)\Big|\lesssim \frac{\sqrt{h}}{\sqrt{q}}+h+O(\min\{\sqrt{q}h,q^\frac 32 h^\frac 32\}).
\end{equation}
On the other hand, 
if $G(p,m_{x_0,q},q)\neq 0$ and $\operatorname{dist}(x_0,\frac{\mathbb Z}{q})=0$, 
we get 
\begin{equation}\label{eqrat}
R_{x_0}\Big(\frac pq+h\Big)-R_{x_0}\Big(\frac pq\Big)
\simeq \frac{\sqrt{h}}{\sqrt{q}} + ih + O(\min\{\sqrt{q}h,q^\frac 32 h^\frac 32\}).
\end{equation}
In particular, if $x_0=\frac PQ\in\mathbb Q$, these conditions are satisfied for $q \in 4Q\mathbb N$, so
\begin{equation}\label{Ho12}
R_{ P/Q} \in \mathcal C^{1/2} \left(p/q\right) \qquad \text{ if } \, \,  q \in 4Q\mathbb N.
\end{equation}

\subsection{Upper and lower bounds for H\"older regularity at irrationals}

To obtain a lower bound for H\"older regularity at an irrational $t$,
we start by recalling that the exponent of irrationality of $t$ is
$$\mu(t)=\sup\{\mu>0 \, : \, |t-\frac pq|\leq \frac{1}{q^\mu} \mbox{ for infinitely many coprime pairs } (p,q)\in \mathbb N\times\mathbb N\}.$$
The approximations by continuous fractions of $t$, denoted by $p_n/q_n$,  satisfy
$$\Big| t-\frac{p_n}{q_n} \Big|=\frac1{q_n^{\mu_n}}
\qquad \text{ and } \qquad 
\mu(t)=\limsup_{n\rightarrow\infty}\mu_n.$$
Since $p_n/q_n \to t$, we have that
\begin{equation}\label{H_Trapped}
\forall h > 0, \qquad \exists n \in \mathbb N \quad : \quad  \Big|t-\frac{p_n}{q_n}\Big|\leq h\leq \Big|t-\frac{p_{n-1}}{q_{n-1}}\Big|.
\end{equation} 
We use \eqref{varrat}, \eqref{H_Trapped} and the property $|t - p_{n-1}/q_{n-1}| \simeq (q_n q_{n-1})^{-1}$ to estimate
\begin{equation*}
\begin{split}
|R_{x_0}(t+h)-R_{x_0}(t)| & 
\leq \Big| R_{x_0}\Big(\frac{p_n}{q_n} + \big(t-\frac{p_n}{q_n}+h\big) \Big)-R_{x_0}\Big(\frac{p_n}{q_n}\Big) \Big| \\
& \qquad \qquad \qquad 
+ \Big| R_{x_0}\Big(\frac{p_n}{q_n}+ \big( t-\frac{p_n}{q_n}\big) \Big) - R_{x_0}\Big(\frac{p_n}{q_n}\Big) \Big| \\
& \lesssim \frac{\sqrt{h}}{\sqrt{q_n}}+h+\min\{\sqrt{q_n}h,q_n^\frac 32h^\frac 32\} \\ 
& \lesssim h^{\frac 12+\frac{1}{2\mu_n}}+h^{\frac 12+\frac{1}{2\mu_{n-1}}} \\
& \lesssim h^{\frac 12+\frac{1}{2\mu}-\delta}, 
\end{split}
\end{equation*}
which holds for all $\delta > 0$. 
Thus,  the H\"older exponent of $R_{x_0}$ at $t$ satisfies
\begin{equation}\label{lowerHo}
\alpha_{R_{x_0}}(t)\geq\frac 12+\frac1{2\mu(t)}, \qquad \forall t \not\in \mathbb Q.
\end{equation}

To obtain an upper bound for H\"older regularity at irrational $t$, we shall fix from now on $x_0=\frac PQ\in\mathbb Q$ and we shall use the poor regularity of $R_{ P/Q}$ at the rationals specified in \eqref{Ho12}, i.e. with denominator in $4Q\mathbb N$.
Indeed, consider the irrationals well approximated by rationals with denominator restricted to $4Q\mathbb N$,
$$ \mathbf{A}_{\mu,Q} = \{ \, t \notin\mathbb Q \, : \,  |t-\frac pq|\leq \frac{1}{q^\mu} \mbox{ for infinitely many coprime pairs } (p,q)\in \mathbb N\times 4Q\mathbb N \, \}.$$
In particular
$$ t\in \mathbf{A}_{\mu,Q} \quad \Longrightarrow \quad  \mu\leq \mu(t)=\sup\{\nu \, : \,  t\in \mathbf{A}_\nu \, \},$$ 
where
$$\mathbf{A}_{\nu} = \{ \, t \notin\mathbb Q \, : \,   |t-\frac pq|\leq \frac{1}{q^\nu} \mbox{ for infinitely many coprime pairs } (p,q)\in \mathbb N\times \mathbb N \, \}.$$
For $t\in \mathbf{A}_{\mu,Q}$, 
we can pick $(p_n,q_n)\in \mathbb N\times 4Q\mathbb N$ and define $\mu_n$ such that
$$\frac{1}{q_n^{\mu_n}} = \Big| t-\frac{p_n}{q_n}\Big|\leq \frac{1}{q_n^\mu}.$$
We set $h_n=t-\frac{p_n}{q_n}$ and we use \eqref{eqrat} to get the lower estimate
$$|R_{x_0}(t+h_n)-R_{x_0}(t)|
= \Big|R_{x_0}\Big(\frac{p_n}{q_n} \Big)-R_{x_0} \Big(\frac{p_n}{q_n}+h_n\Big)\Big| 
\gtrsim \frac{\sqrt{h_n}}{\sqrt{q_n}}=h_n^{\frac 12+\frac1{2\mu_n}}\geq h_n^{\frac 12+\frac1{2\mu}}.$$
Therefore, by combining this with \eqref{lowerHo} we obtain upper and lower bounds for the H\"older regularity at irrationals of $\mathbf{A}_{\mu,Q}$ given by
\begin{equation}\label{lower-upperHo}
 \frac 12+\frac1{2\mu(t)}\leq \alpha_{R_{x_0}}(t) \leq \frac 12+\frac1{2\mu},
 \qquad \forall t\in \mathbf{A}_{\mu,Q}.
\end{equation}

\subsection{Measuring the iso-H\"older sets}
In \eqref{lower-upperHo} the H\"older exponent would be completely determined if $\mu(t)=\mu$.
However, this is not the case for all $t\in \mathbf{A}_{\mu,Q}$. 
To fix this, we shall remove the points where the H\"older exponent is larger than $\mu$ by introducing the sets
$$\mathbf{B}_{\mu,Q} = \mathbf{A}_{\mu,Q} \setminus \Big(\bigcup_{\epsilon>0}{\bf{A_{\mu+\epsilon}}}\Big).$$
From $\mathbf{B}_{\mu,Q} \subset \mathbf{A}_{\mu}\setminus\Big(\bigcup_{\epsilon>0}\mathbf{A}_{\mu + \epsilon} \Big)$ and the definition of $\mu(t)$, 
we get that
$$t\in \mathbf{B}_{\mu,Q} \quad \Longrightarrow \quad \mu(t)=\mu. $$
Therefore, by \eqref{lower-upperHo} and by the combination of \eqref{lowerHo}
with $\mu(t)=\sup\{\, \nu \, : \,  t\in \mathbf{A}_{\nu} \}$ we get
\begin{equation}\label{isoHosets}
\mathbf{B}_{\mu,Q} \subset \Big\{ \, t \, : \, \alpha_{R_{x_0}}(t) = \frac 12+\frac1{2\mu} \, \Big\}\subset \mathbf{A}_{\mu-\epsilon},
\qquad  \forall\epsilon>0.
\end{equation}
Since Theorem~\ref{th} claims that 
$$\dim_{\mathcal H} \Big\{ \, t \, : \, \alpha_{R_{x_0}}(t)=\frac 12+\frac1{2\mu} \, \Big\}=\frac 2\mu,$$
and since  $\dim_{\mathcal H}\mathbf{A}_\nu= \frac 2{\nu}$  for $\nu \geq 2$ by the Jarn\'{i}k-Besicovitch theorem,
it suffices to show
$$\dim_{\mathcal H}\mathbf{B}_{\mu,Q}\geq \frac 2\mu.$$
Since $\mathcal H^\frac 2\mu(\mathbf{A}_{\mu+\frac 1n})=0$ for all $n \in \mathbb N$,
we can write
$$\mathcal H^\frac 2\mu (\mathbf{B}_{\mu,Q}) 
= \mathcal H^\frac 2\mu (\mathbf{A}_{\mu,Q}) - \lim_{n\rightarrow\infty}\mathcal H^\frac 2\mu(\mathbf{A}_{\mu+\frac 1n})
= \mathcal H^\frac 2\mu (\mathbf{A}_{\mu,Q}),$$
so it is enough to prove $\mathcal H^\frac 2\mu (\mathbf{A}_{\mu,Q})>0$. 
We will actually prove  
\begin{equation}\label{H2/mu}
\mathcal H^\frac 2\mu (\mathbf{A}_{\mu,Q})= \infty, 
\end{equation}
which will complete the proof.

First, we use the Duffin-Schaeffer theorem \cite{KoukoulopoulosMaynard2020}, which states that if 
\begin{equation}\label{DS_Condition}
\sum_{q= 1}^\infty \psi(q)\varphi(q)=\infty,
\end{equation}
where $\varphi$ is Euler's totient function\footnote{Euler's totient function is $\varphi(q) = \# \{\, 1 \leq m \leq q \, : \,  \operatorname{gcd}(m,q) = 1 \, \}$.}
and $\psi$ is an arbitrary function, then the set
$$\mathbf{A}_\psi = \Big\{ \,  t \, : \,  \big|t-\frac pq \big|\leq \psi(q) \, \text{ for infinitely many coprime pairs } (p,q) \in \mathbb N\times \mathbb N \, \Big\}$$
has Lebesgue measure $1$. 
The function
 $\psi(q):=\frac{\mathbbm 1_{4Q\mathbb N}(q)}{q^2}$ satisfies the hypothesis \eqref{DS_Condition} and $\mathbf{A}_\psi = \mathbf{A}_{2,Q}$, so the Duffin-Schaeffer theorem implies
\[
\mathcal L(\mathbf{A}_{2,Q})
=1.
\]
Once we know this,
we use the Mass Transference Principle of Beresnevich and Velani \cite{BeresnevichVelani2006}, 
which says that if $B(x_n,r_n)$ is a sequence of balls in $[0,1]^d$ with $r_n\rightarrow 0$,
and if for some $\alpha < d$ we have
$$\mathcal L\Big(\limsup_n B\big(x_n, r_n^{\alpha/d}\big)\Big) = 1,$$
then 
$$ \dim_{\mathcal H} \Big( \limsup_n B(x_n, r_n) \Big) \geq \alpha,
\qquad \text{ and } \qquad  \mathcal H^\alpha \Big(\limsup_n B(x_n, r_n)\Big) = \infty.$$
In our case $d=1$, since we have obtained 
\begin{equation}
\begin{split}
1 = \mathcal L(\mathbf{A}_{2,Q})
& = \mathcal L\Bigg(\limsup_q \bigcup_{p\leq q,(p,q)=1} B\Big(\frac pq, \frac{\mathbbm 1_{4Q\mathbb N}(q)}{q^2}\Big) \, \Bigg) \\
& = \mathcal L \Bigg(  \limsup_q \bigcup_{p\leq q,(p,q)=1} B\Big(\frac pq,\Big(  \frac{\mathbbm 1_{4Q\mathbb N}(q)}{q^\mu} \Big)^{2/\mu}\Big) \,  \Bigg), 
\end{split}
\end{equation}
applying the result above with $\alpha= 2/ \mu$ we get
$$
\mathcal H^\frac 2\mu \Bigg(\limsup_q \bigcup_{p\leq q,(p,q)=1}B\Big(\frac pq, \frac{\mathbbm 1_{4Q\mathbb N}(q)}{q^\mu}\Big)\Bigg)
= \mathcal H^\frac 2\mu(\mathbf{A}_{\mu,Q})
= \infty, 
$$
which is what we wanted to prove in \eqref{H2/mu}. The proof is complete.

\bibliographystyle{plain}

\bibliography{biblio}

\end{document}